\documentclass[%
 reprint,
 amsmath,amssymb,
 aps,
]{revtex4-2}

\usepackage{graphicx}% Include figure files
\usepackage{dcolumn}% Align table columns on decimal point
\usepackage{bm}% bold math
\begin{document}

%\preprint{APS/123-QED}

\title{Bootstrapping the Stochastic Resonance}% Force line breaks with \\

\author{Minjae Cho}
 \email{minjae@princeton.edu}
\affiliation{
 Princeton Center for Theoretical Science, Princeton University,
 \\Princeton, NJ 08544, USA
}

\begin{abstract}
Stochastic resonance is a phenomenon where a noise of appropriate intensity enhances the input signal strength. In this work, by employing the recently developed convex optimization methods in the context of dynamical systems and stochastic processes, we derive rigorous two-sided bounds on the expected power at the input signal frequency for the prototypical example of stochastic resonance, the double-well potential with periodic forcing and Gaussian white noise.
%\begin{description}
%\end{description}
\end{abstract}

\maketitle

%\tableofcontents

\section{\label{sec:introduction}Introduction}
In contrast to the common intuition that the presence of noise would always lead to a decrease in the power of signals of interest, there are numerous systems where the former leads to an amplification of the latter, a phenomenon called stochastic resonance. Since first introduced in \cite{RBenzi_1981,1982Tell...34...10B,1ad8f307-962d-3238-b7eb-99405febec8a,RBenzi_1985} as a potential explanation of the recurring ice ages, stochastic resonance has been studied and demonstrated in various fields of science, ranging over classical and quantum physics, optics, electronics, medical physics, chemical reactions, biology, neuroscience, etc (see \cite{mcdonnell2009} for a survey on the extensive list of fields where stochastic resonance made appearances, and \cite{RevModPhys.70.223,ThomasWellens_2004} for comprehensive reviews on the topic). The detailed mechanisms of stochastic resonance may differ from one specific system to another, but a common ingredient in many cases is the existence of the characteristic rate which is a function of the noise intensity. By tuning the latter, one may obtain the characteristic rate which is similar to (two times) the frequency of the input signal, thus leading to a resonance effect.

To illustrate the idea further, we shall focus on the prototypical example where an overdamped Brownian particle is placed in a double-well potential with periodic forcing, whose motion is described by a random variable $X(t)$ taking values $x(t)\in\mathbb{R}$ and obeying the stochastic differential equation (SDE)
\begin{equation}\label{eq:SRSDE}
    dX=(-V'(X)+A \text{cos}(\Omega t))dt+\sqrt{2D}dW,
\end{equation}
where $V(x)=-\frac{1}{2}x^2+\frac{1}{4}x^4$ and $W$ is the standard Wiener process. The parameters $A,\Omega,$ and $D$ stand for the amplitude and frequency of the periodic forcing term and the noise intensity respectively. When $A=0$ and $D=0$, there are two stable fixed points at $x=\pm1$, and the particle motion is confined to either the left $(x<0)$ or the right well $(x>0)$ depending on its initial position except for the unstable fixed point $x=0$. When $A=0$ but $D\neq0$, random hopping between the left and the right wells happens at the celebrated Kramer's rate $r_K$ \cite{KRAMERS1940284} approximately given by
\begin{equation}
    r_K\approx \frac{1}{\sqrt{2}\pi}\text{exp}\left(-\frac{\Delta V}{D}\right),
\end{equation}
where $\Delta V$ is the height of the potential barrier which in this case equals $\frac{1}{4}$. The resulting motion $X(t)$ is zero when averaged over the noise.

\begin{figure}[b]
\includegraphics[width=0.5\textwidth]{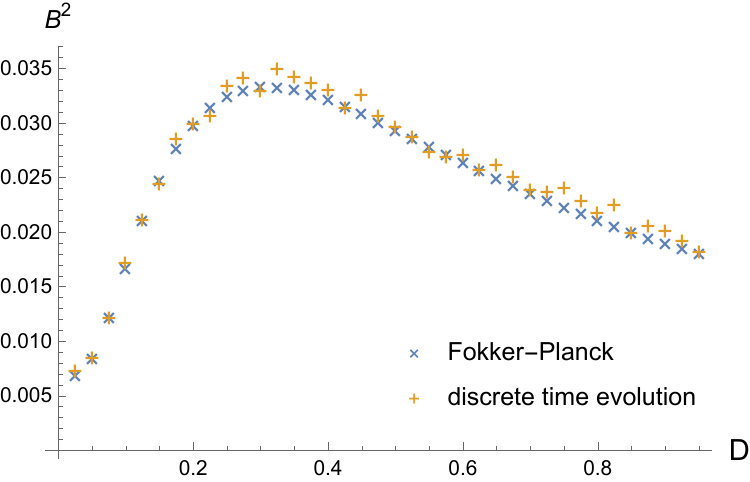}% Here is how to import EPS art
\caption{\label{fig:discreteTime} $B^2$ as a function of the noise intensity $D$, with $A=0.3$ and $\Omega=0.5$. Each point for the discrete time evolution is obtained by averaging over 100 independent discrete time evolutions of (\ref{eq:SRSDE}) using Mathematica's \cite{Mathematica} \texttt{ItoProcess} function with \texttt{Euler-Maruyama} method. Each point for the FP equation is obtained by numerically solving (\ref{eq:FP}) using Mathematica's \texttt{NDSolve} function with \texttt{PDEdiscretization} method with \texttt{MinPoints}=$2\times10^5$ over the domain $(t,x)\in[0,100]\times[-5,5]$.}
\end{figure}

For $0<A<\frac{2}{3\sqrt{3}}$, which is the regime of interest for this work, the periodic forcing term $A\text{cos}(\Omega t)$ is not big enough to allow the particle to travel along both wells in the absence of the noise $(D=0)$. When the noise is turned on $D>0$ in addition, ergodicity implies that after a long time, evolution of any initial $x(t=0)$ asymptotes to a motion whose noise-averaged expectation value $\langle X(t)\rangle$ is periodic in $t$ with periodicity $T_\Omega=\frac{2\pi}{\Omega}$ \cite{P.Jung_1989,PhysRevA.44.8032}, as a result of the Floquet theorem. Therefore, $\langle X(t)\rangle$ allows for the Fourier series expansion
\begin{equation}\label{eq:Fourier}
    \langle X(t)\rangle = 2\sum_{n=0}^\infty \left(a_n \text{cos}(n\Omega t)+b_n \text{sin}(n\Omega t) \right).
\end{equation}
The amplitude for the first frequency mode, $B=\sqrt{a_1^2+b_1^2}$, is one of the measures of stochastic resonance, and $B^2$ is proportional to the power spectral density evaluated at the frequency $\Omega$. It is a function of the noise intensity $D$ which may obtain a maximum at a nonzero finite value of $D$ where the two characteristic rates $r_\Omega=\frac{\Omega}{\pi}$ and $r_K$ become similar (but not necessarily the same). This phenomenon of signal enhancement is the hallmark of stochastic resonance, which is illustrated in Fig.~\ref{fig:discreteTime} for the example discussed here.

Traditionally, there have been two main methods to perform an explicit numerical computation of quantities such as $B$. The first is the discrete time evolution of SDE where at each time step, a random number from the noise term is drawn. In order to perform the noise-averaging and thus compute $\langle X(t)\rangle$, the time evolution is repeated many times independently and then averaged. Due to the inherent stochastic nature of SDE, the statistical variance of the sampled paths cannot be reduced to an arbitrarily small number. In contrast, the second method based on the Fokker-Planck (FP) equation is free from having to draw a random number from the noise term. It describes the evolution of the probability density $p(x,t)$ which provides the expectation values under the noise-averaging. The FP equation corresponding to the SDE (\ref{eq:SRSDE}) is given by
\begin{equation}\label{eq:FP}
    \partial_t p(x,t) = -\partial_x \left(\left(-V'(x)+A \text{cos}(\Omega t)\right)p(x,t)\right)+D\partial_x^2p(x,t).
\end{equation}
By numerically solving (\ref{eq:FP}) for large enough time domain so that $p(x,t)$ approaches an equilibrium density $p_{eq}(x,t)$, one may rewrite (\ref{eq:Fourier}) as
\begin{subequations}
\label{eq:expEq}
\begin{equation}
\langle X(t)\rangle = \int_{-\infty}^\infty dx p_{eq}(x,t) x,
\end{equation}
\begin{equation}
a_1 = \frac{1}{T_\Omega}\int_T^{T+T_\Omega} dt\text{cos}(\Omega t)\int_{-\infty}^\infty dx p_{eq}(x,t) x,
\end{equation}
\begin{equation}
b_1 = \frac{1}{T_\Omega}\int_T^{T+T_\Omega} dt\text{sin}(\Omega t)\int_{-\infty}^\infty dx p_{eq}(x,t) x,
\end{equation}
\end{subequations}
for any $T\in\mathbb{R}$. Therefore, the FP equation provides a direct access to the equilibrium probability density from which noise-averaged quantities may be obtained. Numerical results obtained from the two aforementioned methods are presented in Fig.~\ref{fig:discreteTime}.

Recently, a new approach to the invariant measures of dynamical and stochastic systems using the convex optimization was introduced in \cite{2015arXiv151205599F}. It was further explored in more examples of SDEs such as nonlinear drift or diffusion in \cite{2018arXiv180708956K} and a minimal model of turbulence in \cite{PhysRevE.107.054114}. The essential idea of the approach is to make use of two properties of the invariant measure, the first being the positivity saying that the measure should be non-negative, and the second being the invariance saying that the measure should be invariant under the time evolution. Imposing such consistency conditions of the objective to derive nontrivial consequences about the physical observables of interest is in general called the bootstrap approach, which has had great successes in several important problems in theoretical physics such as conformal bootstrap \cite{Rattazzi:2008pe,El-Showk:2012cjh,Kos:2016ysd,Poland:2018epd}.

Aforementioned two properties can be combined into a semidefinite programming (SDP) problem where one minimizes/maximizes an objective expectation value of a function of random variables over all the invariant measures. In the case of stochastic resonance where the SDE is not autonomous, there is no such invariant measure. However, as non-autonomous SDEs can be made autonomous by introducing more variables, one may consider the higher-dimensional autonomous SDEs which possess an invariant measure, whose expectation values agree with those evaluated using the equilibrium probability density $p_{eq}(x,t)$ upon averaging over the time interval $t\in(T,T+T_\Omega]$, as explained for example in \cite{P.Jung_1989,PhysRevA.44.8032}.

In this work, we employ the convex optimization approach to study such expectation values, which in particular include $a_1$ and $b_1$. Because it relies on the properties that $p_{eq}(x,t)$ must satisfy, the minima/maxima obtained from the optimization problem provide rigorous lower/upper bounds on the expectation values, which may be very close to each other in some cases. For stochastic resonance, such rigorous bounds on the signal amplification (which is proportional to $B$) as a function of the noise intensity is both conceptually and quantitatively desirable because lower bounds imply how much amplification is guaranteed around the noise intensity where resonance takes place, while upper bounds provide its limitations.

\section{SDP FORMULATION OF SDE}
\subsection{Autonomous case}
We first review how to formulate the SDP problem for the stationary (which is a synonym for the invariant) probability measure of an autonomous SDE, as outlined in \cite{2015arXiv151205599F,2018arXiv180708956K}. Given an autonomous SDE of the form
\begin{equation}\label{eq:genSDE}
    d\mathbf{Y}=\mu(\mathbf{Y})dt+\sigma(\mathbf{Y})d\mathbf{W},
\end{equation}
where $\mathbf{Y}$ is a $n$-dimensional vector-valued random variable taking values $\mathbf{y}(t)\in\mathbb{R}^n$, $\mu(\mathbf{Y})$ is a $n$-dimensional vector drift term, $\sigma(\mathbf{Y})$ is a $n\times k$ matrix diffusion term, and $\mathbf{W}$ is a $k$-dimensional vector-valued standard Wiener process, the corresponding FP equation for the probability density $p(\mathbf{y},t)$ is given by
\begin{equation}
    \partial_t p(\mathbf{y},t)=-\sum_{i=1}^n\partial_{\mathbf{y}_i}\left(\mu_i(\mathbf{y})p(\mathbf{y},t)\right)+\sum_{i,j=1}^n\partial_{\mathbf{y}_i}\partial_{\mathbf{y}_j}\left(D_{ij}(\mathbf{y})p(\mathbf{y},t)\right),
\end{equation}
where $D_{ij}(\mathbf{y})=\frac{1}{2}\left(\sigma(\mathbf{y})\sigma(\mathbf{y})^T\right)_{ij}$. A stationary probability density $p_s(\mathbf{y})$ then solves
\begin{equation}\label{eq:sFP}
    -\sum_{i=1}^n\partial_{\mathbf{y}_i}\left(\mu_i(\mathbf{y})p_s(\mathbf{y})\right)+\sum_{i,j=1}^n\partial_{\mathbf{y}_i}\partial_{\mathbf{y}_j}\left(D_{ij}(\mathbf{y})p_s(\mathbf{y})\right)=0.
\end{equation}
and we denote the corresponding expectation value of a function $f(\mathbf{Y})$ by $\langle\langle f(\mathbf{Y}) \rangle\rangle = \int d\mathbf{y} p_s(\mathbf{y})f(\mathbf{y})$ for all $f(\mathbf{Y})$ such that the integral makes sense. Using integration by parts, (\ref{eq:sFP}) implies
\begin{equation}\label{eq:sFPex}
    \bigg\langle\bigg\langle \sum_{i=1}^n\mu_i(\mathbf{Y})\partial_{\mathbf{Y}_i}f(\mathbf{Y})+\sum_{i,j=1}^nD_{i,j}(\mathbf{Y})\partial_{\mathbf{Y}_i}\partial_{\mathbf{Y}_j}f(\mathbf{Y}) \bigg\rangle\bigg\rangle=0.
\end{equation}

We denote the monomials in $\mathbf{Y}$ by $\mathbf{Y}^{\alpha}=\prod_{i=1}^n\mathbf{Y}_i^{{\alpha}_i}$ and define $|{\alpha}|=\sum_{i=1}^n{\alpha}_i$. Also, we denote the set of monomials up to degree $d$ by $M_d=\{\mathbf{Y}^{\alpha}|~|\alpha|\leq d\}$. For simplicity, we assume that $\mu(\mathbf{Y})$ is a polynomial of degree $d_\mu$ and $\sigma(\mathbf{Y})$ is a polynomial of degree $d_\sigma$. Positivity of the stationary probability density $p_s(\mathbf{y})$ implies
\begin{equation}
    {\cal M}^d\succeq0,~~~\forall d\in\mathbb{N},
\end{equation}
where the matrix ${\cal M}^d$ is defined by its matrix elements $\left({\cal M}^d\right)_{\alpha,\beta}\equiv \big\langle\big\langle \mathbf{Y}^\alpha\mathbf{Y}^\beta  \big\rangle\big\rangle$ with $|\alpha|\leq d$ and $|\beta|\leq d$. If the objective we hope to minimize is given by $\sum_{\alpha}c_\alpha\langle\langle\mathbf{Y}^\alpha\rangle\rangle$ for a given real coefficients $c_\alpha$, the SDP problem at degree $d$ is given by
\begin{eqnarray}\label{eq:SDPst}
    &&\text{minimize }\sum_{\alpha}c_\alpha\langle\langle\mathbf{Y}^\alpha\rangle\rangle,~\nonumber
    \\
    &&\text{subject to } {\cal M}^d\succeq0~\text{and}\nonumber
    \\
    &&\bigg\langle\bigg\langle \sum_{i=1}^n\mu_i(\mathbf{Y})\partial_{\mathbf{Y}_i}\mathbf{Y}^\alpha+\sum_{i,j=1}^nD_{i,j}(\mathbf{Y})\partial_{\mathbf{Y}_i}\partial_{\mathbf{Y}_j}\mathbf{Y}^\alpha \bigg\rangle\bigg\rangle=0,\nonumber
    \\
    &&\forall \alpha~\text{such that}~|\alpha|+\text{max}(d_\mu-1,2d_\sigma-2)\leq 2d.
\end{eqnarray}
The resulting minimum value of the objective provides a rigorous lower bound on $\sum_{\alpha}c_\alpha\langle\langle\mathbf{Y}^\alpha\rangle\rangle$ realized by all stationary probability densities. An upper bound can be obtained similarly. The convergence of the bounds to the actually existing expectation values as $d$ increases can be shown in case where the support of $p_s(\mathbf{y})$ is compact and the corresponding additional positivity constraints are added to (\ref{eq:SDPst}), as discussed in \cite{2018PhLA..382..382T,2018arXiv180708956K}. To be precise, in order to obtain absolutely rigorous bounds, one should further employ interval arithmatic, which we do not attempt in this work.

\subsection{Non-autonomous to autonomous}
The SDE (\ref{eq:SRSDE}) is non-autonomous due to the presence of periodic forcing term. By introducing non-random variables $y(t)=\text{cos}(\Omega t)$ and $z(t)=\text{sin}(\Omega t)$, (\ref{eq:SRSDE}) can be rewritten as
\begin{eqnarray}\label{eq:SRauto}
    &&dX = (-V'(X)+Ay)dt+\sqrt{2D}dW,\nonumber
    \\
    &&dy = -\Omega zdt,~~~dz = \Omega ydt,~~y(t)^2+z(t)^2=1.
\end{eqnarray}
This is of the form (\ref{eq:genSDE}) with $\mathbf{Y}=(X,y,z)^T$, $\mu(\mathbf{Y})=(-V'(X)+Ay,-\Omega z,\Omega y)^T$, $\sigma(\mathbf{Y})=(\sqrt{2D},0,0)^T$, and $k=1$, but also with the extra constraint $y(t)^2+z(t)^2=1$. There is a $\mathbb{Z}_2$-symmetry which acts as $(X,y,z)\rightarrow-(X,y,z)$. Due to the uniqueness of the stationary probability density $p_s(\mathbf{y})$ of the system \cite{MATTINGLY2002185}, expectation values of $\mathbb{Z}_2$-odd variables vanish. Therefore, SDP formulation of (\ref{eq:SRauto}) is given by
\begin{eqnarray}\label{eq:SDPSR}
    &&\text{minimize }\sum_{\alpha}c_\alpha\langle\langle\mathbf{Y}^\alpha\rangle\rangle,~\nonumber
    \\
    &&\text{subject to } {\cal M}^d\succeq0,\nonumber
    \\
    &&\big\langle\big\langle \mathbf{Y}^\beta \left(y^2+z^2-1\right) \big\rangle\big\rangle=0,~\forall \beta~\text{s.t.}~|\beta|\leq2d-2,\nonumber
    \\
    &&\big\langle\big\langle \mathbf{Y}^\gamma \big\rangle\big\rangle=0,~\forall \gamma~\text{s.t.}~|\gamma|~\text{odd and }|\gamma|\leq2d,~\text{and}\nonumber
    \\
    &&\bigg\langle\bigg\langle \sum_{i=1}^n\mu_i(\mathbf{Y})\partial_{\mathbf{Y}_i}\mathbf{Y}^\alpha+\sum_{i,j=1}^nD_{i,j}(\mathbf{Y})\partial_{\mathbf{Y}_i}\partial_{\mathbf{Y}_j}\mathbf{Y}^\alpha \bigg\rangle\bigg\rangle=0,\nonumber
    \\
    &&\forall \alpha~\text{s.t.}~|\alpha|+\text{max}(d_\mu-1,2d_\sigma-2)\leq 2d.
\end{eqnarray}
The condition $\big\langle\big\langle \mathbf{Y}^\beta \left(y^2+z^2-1\right) \big\rangle\big\rangle=0$ enforces that the support of the stationary probability density $p_s(\mathbf{y})$ is localized to $y^2+z^2=1$. A similar SDP problem where non-autonomous deterministic systems are transformed into autonomous deterministic systems was recently discussed in \cite{Doering2020OptimalTA}. The relation between the expectation values under $p_s(\mathbf{y})$ and those under $p_{eq}(x,t)$ are given by \cite{P.Jung_1989,PhysRevA.44.8032}
\begin{equation}
    \big\langle\big\langle f(X,y,z) \big\rangle\big\rangle=\frac{1}{T_\Omega}\int_{T}^{T+T_\Omega}dt \big\langle f(X(t),\text{cos}(\Omega t),\text{sin}(\Omega t)) \big\rangle,
\end{equation}
for any $T\in\mathbb{R}$ and all $f(X,y,z)$ where both sides are well-defined. In particular, we will be interested in the following quantities
\begin{eqnarray}
    &&P=\frac{1}{T_\Omega}\int_T^{T+T_\Omega}dt\big\langle X(t)^2\big\rangle=\big\langle\big\langle X^2 \big\rangle\big\rangle,\nonumber
    \\
    &&a_1=\big\langle\big\langle Xy \big\rangle\big\rangle,~~b_1=\big\langle\big\langle Xz \big\rangle\big\rangle,\nonumber
    \\
    &&B^2=a_1^2+b_1^2=\big\langle\big\langle Xy \big\rangle\big\rangle^2+\big\langle\big\langle Xz \big\rangle\big\rangle^2.
\end{eqnarray}
$P$ is the time average of the expected power. Therefore, another measure of the stochastic resonance is given by $R\equiv \frac{B^2}{P}$.

\begin{figure*}
\centering
\includegraphics[width=0.765\textwidth]{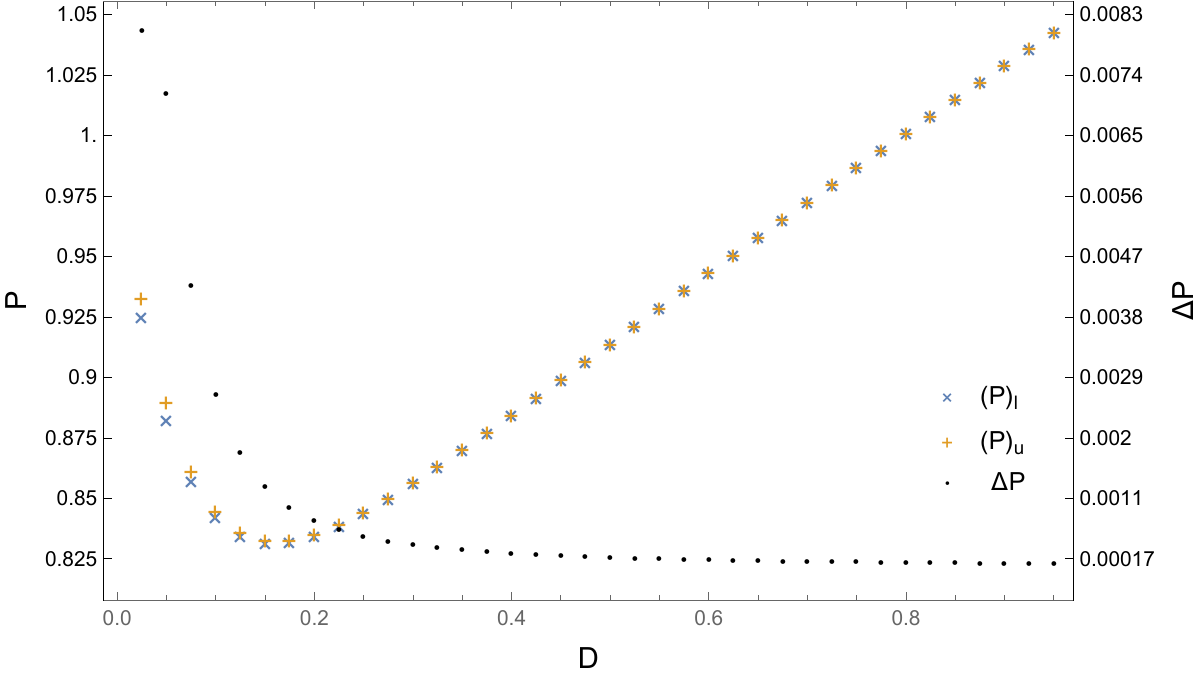}\vspace{0.5cm}
\includegraphics[width=0.75\textwidth]{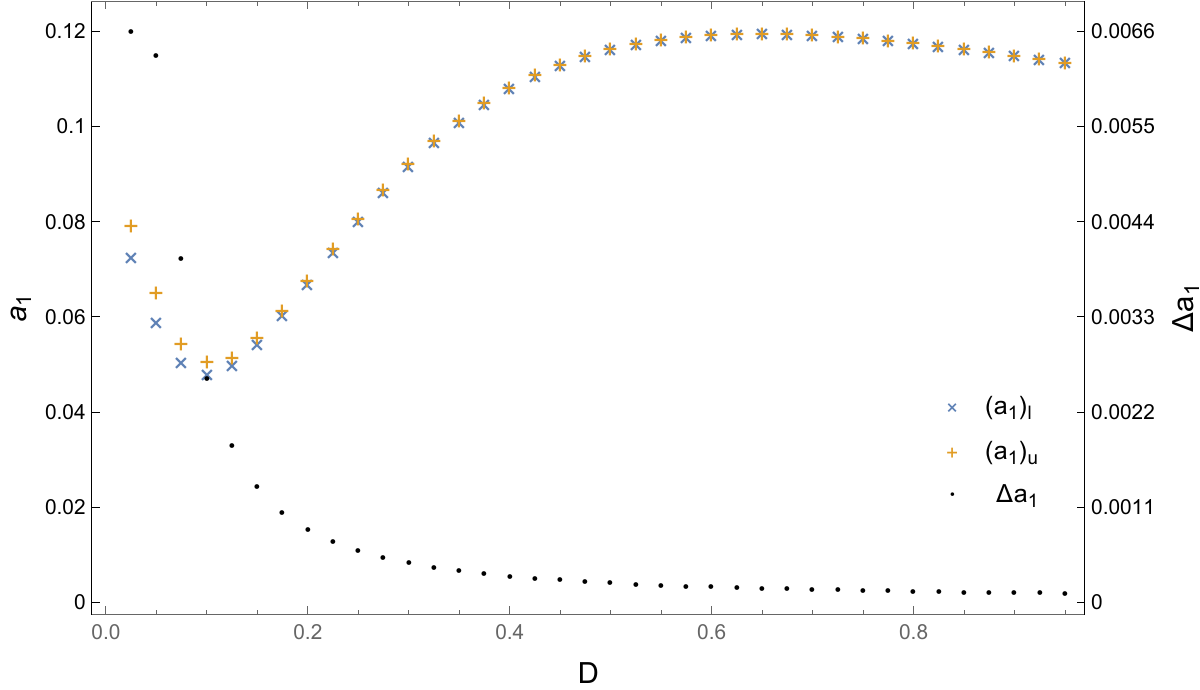}\vspace{0.5cm}
\includegraphics[width=0.75\textwidth]{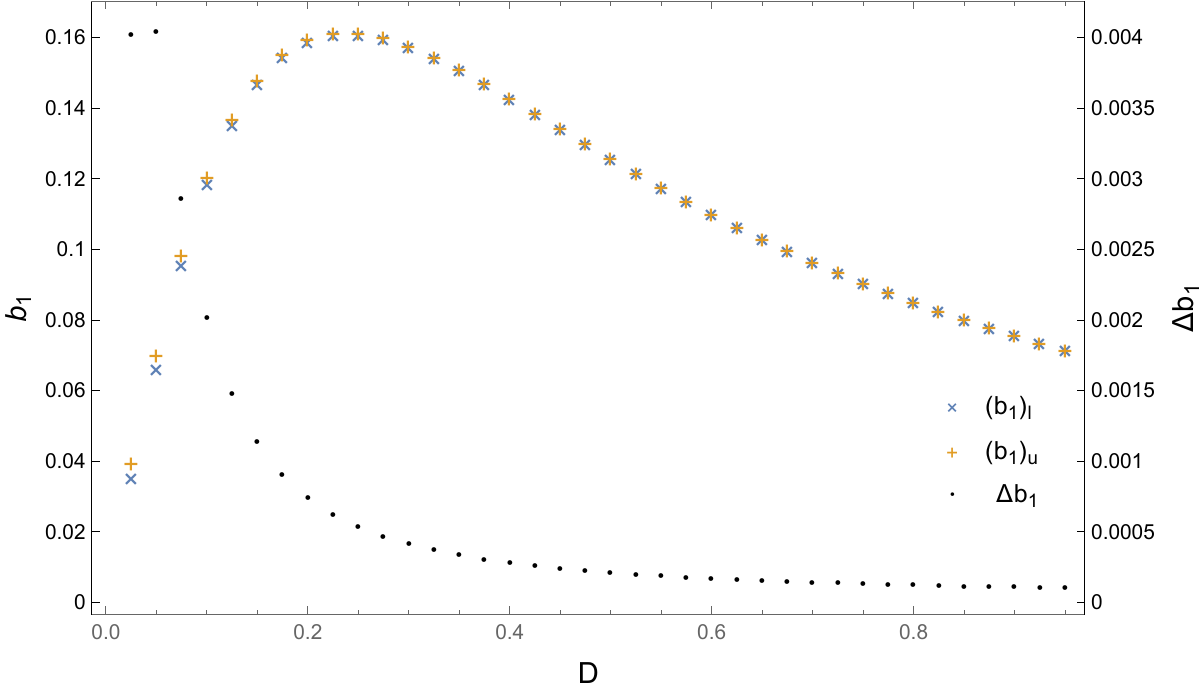}
\caption{\label{fig:deg12} Results of SDP (\ref{eq:SDPSR}) for $P,a_1,$ and $b_1$ at $d=12$. Lower and upper bounds, and their differences are displayed for different values of the noise intensity $D$.}
\end{figure*}

\begin{figure}
\centering
\vspace{0.5cm}
\includegraphics[width=0.45\textwidth]{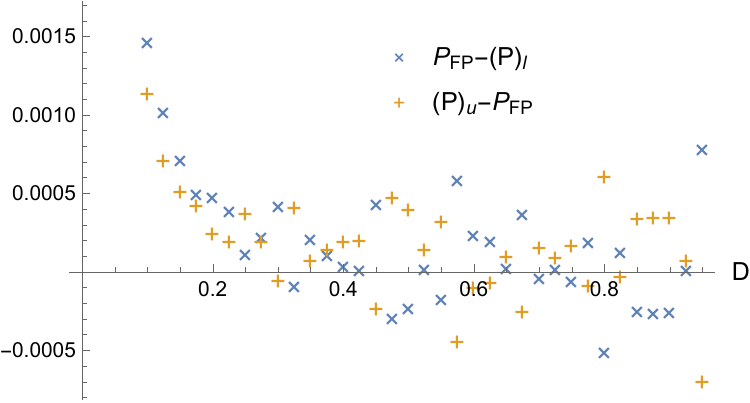}\vspace{0.5cm}
\includegraphics[width=0.45\textwidth]{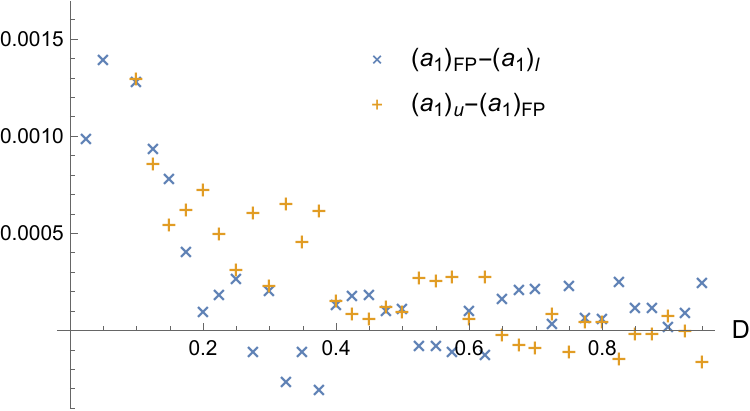}\vspace{0.5cm}
\includegraphics[width=0.45\textwidth]{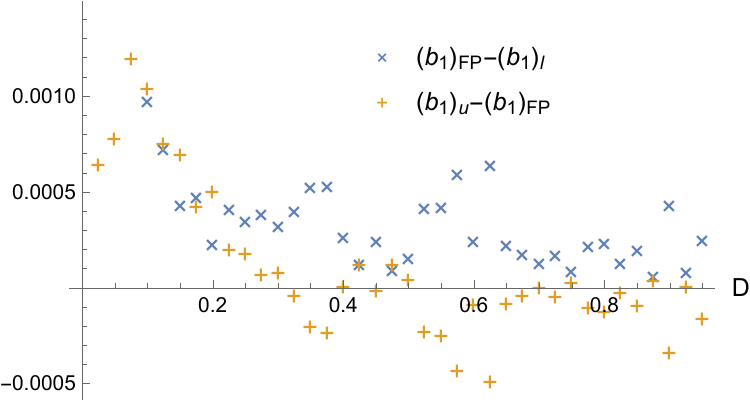}
\caption{\label{fig:againstFP} Comparisons of bounds obtained at $d=12$ and numerical values obtained by the FP equation (denoted with the subscript $FP$). When the point displayed takes a negative value, the corresponding FP value lies outside the range of values allowed by SDP bounds. Points not displayed for small $D$ are positive and bigger than the scales shown in the plots.}
\end{figure}

\section{RESULTS}
Given $d\in\mathbb{N}$, the SDP (\ref{eq:SDPSR}) and the analogous maximization problem, where the objective is taken to be one of $P,a_1$, and $b_1$, provide lower (labeled by subscript $l$) and upper (labeled by subscript $u$) bounds on each of the objectives: $\left(P\right)_l\leq P\leq \left(P\right)_u$, $\left(a_1\right)_l\leq a_1 \leq \left(a_1\right)_u$, and $\left(b_1\right)_l\leq b_1 \leq \left(b_1\right)_u$. When $\left(a_1\right)_l\geq0$ and $\left(b_1\right)_l\geq0$, they further provide bounds on $B^2$ and $R$ by $\left(B^2\right)_l=\left(a_1\right)_l^2+\left(b_1\right)_l^2\leq B^2\leq\left(a_1\right)_u^2+\left(b_1\right)_u^2=\left(B^2\right)_u $ and $\left(R\right)_l=\frac{\left(B^2\right)_l}{\left(P\right)_u} \leq R\leq\frac{\left(B^2\right)_u}{\left(P\right)_l}=\left(R\right)_u $. For each of $P,a_1,b_1,B^2,$ and $R$, we denote the difference between upper and lower bounds by $\Delta P, \Delta a_1,\Delta b_1,\Delta B^2,$ and $\Delta R$ respectively. We fix $A=0.3$ and $\Omega=0.5$.

SDP (\ref{eq:SDPSR}) can be solved using standard double-precision solvers such as MOSEK \cite{mosek} up to $d\lesssim9$, but a higher-precision solver is needed for higher values of $d$. We used the double-double precision solver SDPA-DD \cite{5612693,10.1063/1.2911696,sdpa} on a laptop, where the solver run time at $d=12$ for a single SDP problem was $\sim6$ minutes \footnote{For $d=12,13,$ and $14$, the parameters of SDPA-DD were taken to be \texttt{epsilonStar}=$10^{-16}$, \texttt{lambdaStar}=$10^3$, \texttt{omegaStar}=$2$, \texttt{lowerBound}=$-10^{10}$, \texttt{upperBound}=$10^{10}$, \texttt{betaStar}=$0.1$, \texttt{betaBar}=$0.2$, \texttt{gammaStar}=$0.8$, and \texttt{epsilonDash}=$10^{-16}$, where the notations follow SDPA manual \cite{sdpaManual}. For $d=15$, all the parameters were the same as above except for \texttt{lambdaStar}=$10^5$, \texttt{omegaStar}=$1.5$, and \texttt{gammaStar}=$0.6$.}. In Fig. \ref{fig:deg12}, lower and upper bounds together with their differences for $P, a_1,$ and $b_1$ at $d=12$ are presented as a function of the noise intensity $D$. The difference between lower and upper bounds decreases rather quickly as $D$ increases, to a degree where it is hardly possible to distinguish them by bare eyes. It seems to suggest that bootstrap approach is more powerful for more noisy case where the time required to reach the equilibrium is also expected to be shorter.

The bounds may be compared with the value obtained by the numerical solution to the FP equation presented in Fig. \ref{fig:discreteTime}. Denoting the FP values with the subscript $FP$, the differences $(v)_{FP}-(v)_l$ and $(v)_u-(v)_{FP}$ for $v\in\{P,a_1,b_1\}$ are presented in Fig. \ref{fig:againstFP}. Points which take negative values indicate the cases where FP value lies outside the region allowed by $d=12$ SDP bounds. It should be noted though that there are more advanced numerical methods to solve FP equations (see for example \cite{Kumar2006,PhysRevE.99.032117}) which may produce more precise results. Nonetheless, the comparative advantage of the bootstrap approach is that the bounds are rigorous.

\begin{figure}
\centering
\includegraphics[width=0.52\textwidth]{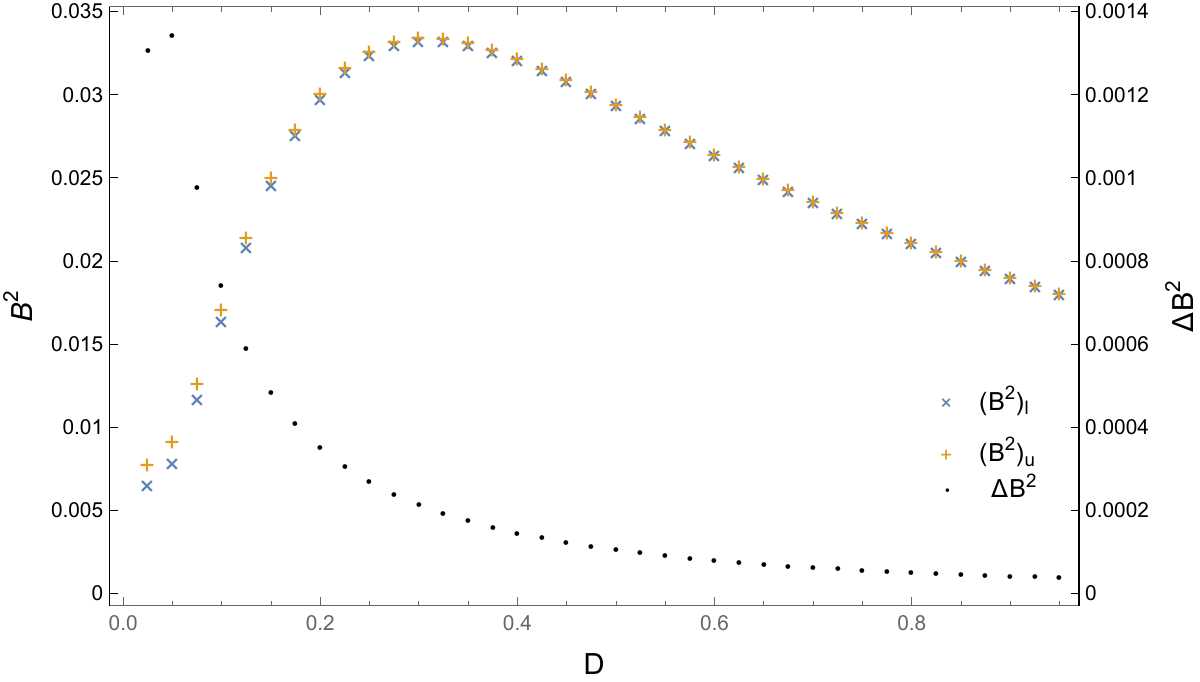}\vspace{0.5cm}
\includegraphics[width=0.52\textwidth]{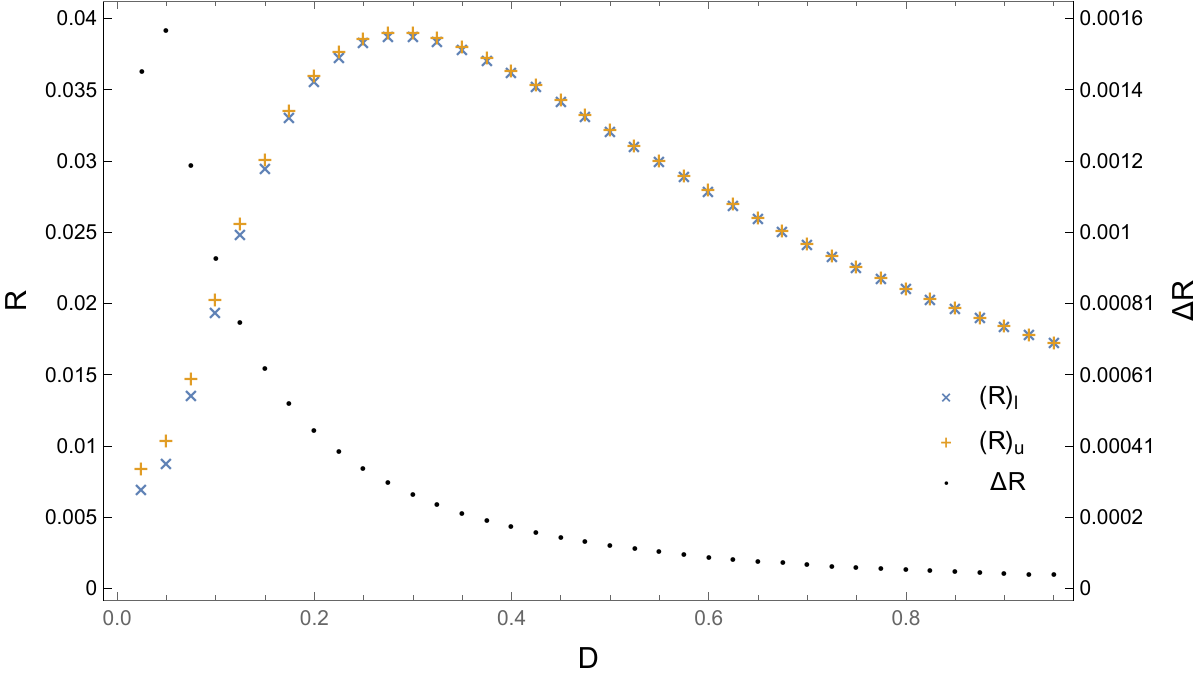}
\caption{\label{fig:B2R} Bounds on $B^2$ and $R$ implied by bounds presented in Fig. \ref{fig:deg12}. Lower and upper bounds, and their differences are displayed for different values of the noise intensity $D$.}
\end{figure}

Using the bounds on $P, a_1,$ and $b_1$, we further obtain bounds on the measures of stochastic resonance $B^2$ and $R$, which are presented in Fig. \ref{fig:B2R}. Stochastic resonance manifests itself as a peak in these measures around $D\sim0.3$.

\begin{figure}
\centering
\includegraphics[width=0.45\textwidth]{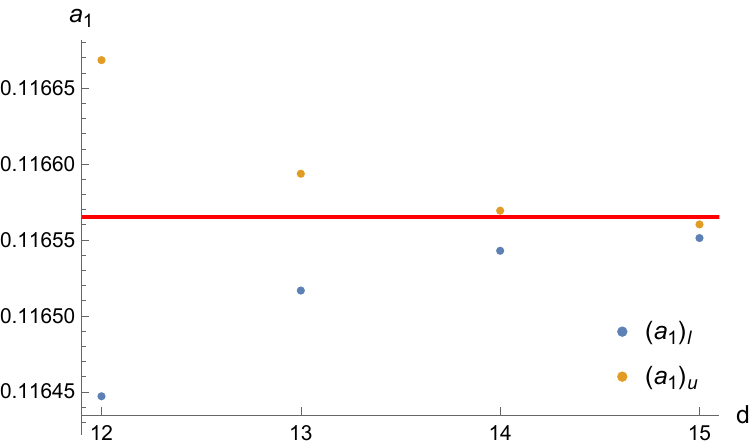}\vspace{0.5cm}
\includegraphics[width=0.45\textwidth]{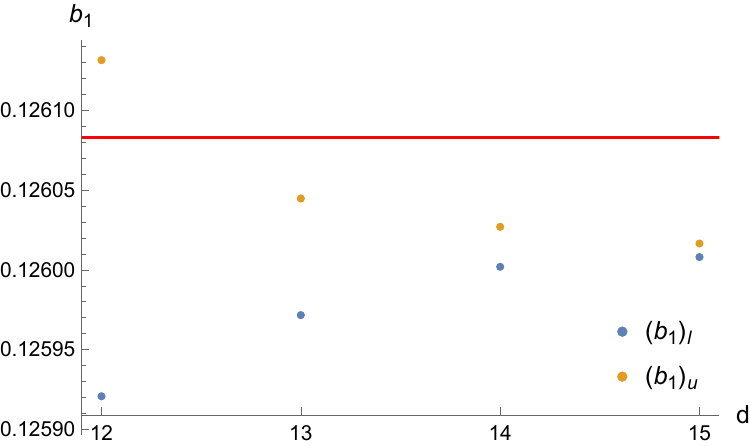}
\caption{\label{fig:higherd} Lower and upper bounds on $a_1$ and $b_1$ for $d=12,13,14$, and $15$ at $D=0.5$, together with the FP values (red solid lines).}
\end{figure}

One may also obtain higher $d$ results, which are presented in Fig. \ref{fig:higherd} for $D=0.5$. SDPA-DD solver run time for $d=13,14$ and $15$ was $\sim13,30,$ and $100$ minutes per point on a laptop respectively. Bounds for $d=15$ up to 8 significant digits where lower bound is rounded down and upper bound is rounded up are given by
\begin{eqnarray}{}
    0.11655120~ \leq ~&&a_1 ~\leq~0.11656018,\nonumber
    \\
    0.12600800~ \leq ~&&b_1 ~\leq~0.12601651.
\end{eqnarray}
The difference between lower and upper bounds in both cases is smaller than $10^{-5}$.

\section{Conclusion}
In this work, we obtained rigorous and tight bounds on the signal amplification for the prototypical example of stochastic resonance using convex optimization. As already illustrated in several works in the literature in recent years, such bootstrap approach to SDEs and random processes provides a direct access to equilibrium probability densities where the rigor is guaranteed and precision may also be present. In particular when the support of the probability density is compact, the usual arguments on the convergence of convex optimization apply and the bounds are guaranteed to converge to physically realized expectation values. Obtaining even higher precision results for the particular example presented in this work seems very doable with the cluster-computing, given that all the results presented in this work were obtained merely on a laptop.

An aspect which in general still requires more studies is the rate of the convergence of such convex optimization methods. In the example discussed in this work, this amounts to asking how fast SDP bounds converge as $d$ increases, as a function of the parameters such as $D, A$ and $\Omega$. Even though the details are not obvious at the moment, it seems reasonable to expect that when the parameters are such that the time it takes for a generic initial distribution to reach the equilibrium distribution is shorter, SDP bounds converge faster. Such expectation may be applicable to convex optimization approach to more general stochastic processes.

\begin{acknowledgments}
MC is supported by the Sam B. Treiman fellowship at the Princeton Center for Theoretical Science. We are grateful to Nicole Shibley for introducing us to stochastic climate models, which led us to learn about stochastic resonance. We are also grateful to Barak Gabai for explaining possible issues of SDPA to us.
\end{acknowledgments}

\nocite{*}

\bibliography{apssamp}% Produces the bibliography via BibTeX.

\end{document}